\documentclass[graybox]{svmult}

\usepackage{amssymb}

\usepackage{mathptmx}       
\usepackage{helvet}         
\usepackage{courier}        
\usepackage{type1cm}        
\usepackage{makeidx}         
\usepackage{graphicx}        
\usepackage{multicol}        
\usepackage[bottom]{footmisc}

\begin{document}

\title*{On the correspondence between mirror-twisted sectors for N=2 supersymmetric vertex operator superalgebras of the form $V \otimes V$ and N=1 Ramond sectors of $V$}
\titlerunning{Mirror-twisted sectors for N=2 VOSAs $V \otimes V$ from N=1 Ramond sectors for $V$}
\author{Katrina Barron}
\institute{Katrina Barron\at University of Notre Dame, Notre Dame, Indiana, USA, \email{kbarron@nd.edu}}
\maketitle

\abstract{Using recent results of the author along with Vander Werf, we present the classification and construction of mirror-twisted modules for N=2 supersymmetric vertex operator superalgebras of the form $V \otimes V$ for the signed transposition mirror map automorphism.  In particular, we show that the category of such mirror-twisted sectors for $V\otimes V$ is isomorphic to the category of N=1 Ramond sectors for $V$.   
}

\section{Introduction}
\label{intro}
In \cite{B-n2twisted},  \cite{B-varna2011}, the author studied twisted modules for N=2 supersymmetric vertex operator superalgebras (N=2 VOSAs) for finite order VOSA automorphisms arising from automorphisms of the N=2 Neveu-Schwarz algebra of N=2 infinitesimal superconformal transformations.  Among such automorphisms is the mirror map.  In \cite{B-n2twisted}, mirror maps were given for N=2 VOSAs of the form $V \otimes V$ where $V$ is an N=1 supersymmetric VOSA of the form $V_L \otimes V_{fer}$, where $V_L$ is a rank $d$ lattice VOSA or the $d$ free boson vertex operator algebra and $V_{fer}$ is the $d$ free fermion VOSA.  In particular, we showed that one of the mirror maps for such an N=2 VOSA, $V \otimes V$, is given by the signed transposition map
\begin{equation}\label{mirror-map} \tilde{\kappa} = (1 \; 2)\ : \ V\otimes V  \longrightarrow   V \otimes V, \qquad \quad   u \otimes v \ \mapsto \  (-1)^{|u||v|} v \otimes u
\end{equation} 
where $|v| = j \, \mathrm{mod} \, 2$ for $v \in V^{(j)}$, with the $\mathbb{Z}_2$-grading of $V$ given by $V = V^{(0)} \oplus V^{(1)}$.

In \cite{B-odd-superpermutation} and \cite{BV-even-superpermutation}, the author along with Vander Werf constructed and classified the cyclic permutation-twisted $V^{\otimes k}$-modules, where $V$ is any VOSA and $k$ is a positive integer.   For $k$ even, this classification is in terms of parity-twisted $V$-modules where the parity automorphism of a VOSA is the map
\begin{equation}\label{parity-map} \sigma \ : \ V  \longrightarrow   V, \qquad \quad   v \ \mapsto \  (-1)^{|v|} v.
\end{equation} 

In this note, we apply the results of \cite{BV-even-superpermutation} to the setting of the mirror map (\ref{mirror-map}) acting on an N=2 supersymmetric VOSA of the form $V \otimes V$, to show that the category of $\tilde{\kappa}$-twisted $(V \otimes V)$-modules is isomorphic to the category of $\sigma$-twisted $V$-modules, which are the N=1 Ramond sectors for the N=1 supersymmetric VOSA, $V$.  This classification also provides an explicit construction of these modules.  

In particular, our result shows that if a representation $M_\sigma$ of the N=1 Ramond algebra is also a parity-twisted modules for a VOSA $V$, where $V \otimes V$ is N=2 supersymmetric, then $M_\sigma$ is also naturally a representation of the mirror-twisted N=2 Neveu-Schwarz algebra.  These results can be used to calculate the graded dimensions for one module in terms of the graded dimensions for the other as shown in Corollary \ref{graded-dimension-corollary} below.   Note that for our results, we do not need to make any assumptions about, for instance, the values of the central charge, the complete reducibility of the representations, or the rationality of the VOSAs.

Certain representations of the N=1 Ramond algebra and related VOSA constructions have previously been studied in, e.g.,  \cite{FQS1985}, \cite{GKO1986}, \cite{FFR}, \cite{Li-twisted}, \cite{S}, \cite{IK1999+},  \cite{Milas}, \cite{Milas-Adamovic2008}.  Certain representations of the mirror-twisted N=2 Neveu-Schwarz algebra have previously been studied in, e.g., \cite{BFK1986}, \cite{Dobrev1987}, \cite{Matsuo1987}, \cite{Kiritsis1988}, \cite{DG2001}, \cite{Gato-Rivera2002}, \cite{IK2004}, \cite{LSZ2010}. 
In particular, the relationship between characters of certain modules for the N=1 Ramond algebra and certain modules for the mirror-twisted N=2 Neveu-Schwarz algebra had previously been observed.  Our explicit isomorphism between mirror-twisted sectors for $V \otimes V$ and  N=1 Ramond sectors for $V$, gives a constructive and overarching explanation of this phenomenon through the theory of VOSAs.

\section{The notions of VOSA and twisted module}
\label{definitions-section}

Following the notation of \cite{B-odd-superpermutation}, \cite{BV-even-superpermutation}, we recall the notion of VOSA and the notions of weak, weak admissible and ordinary $g$-twisted $V$-module for a VOSA, $V$, and an automorphism $g$ of $V$ of finite order.

Let $x, x_0, x_1, x_2,$ denote commuting independent formal variables.
Let $\delta (x) = \sum_{n \in \mathbb{Z}} x^n$.  Expressions such as $(x_1 -
x_2)^n$ for $n \in \mathbb{C}$ are to be understood as formal power series expansions in
nonnegative integral powers of the second variable.

\begin{definition}
A {\it vertex operator superalgebra} is a $\frac{1}{2}\mathbb{Z}$-graded (by weight) vector space $
V=\bigoplus_{n\in \frac{1}{2} \mathbb{Z}}V_n$, 
satisfying ${\rm dim} \, V_n < \infty$ and $V_n = 0$ for $n$ sufficiently negative, 
that is also $\mathbb{Z}_2$-graded by {\it sign}, $V = V^{(0)} \oplus V^{(1)}$, with $V^{(j)} = \bigoplus_{n \in \mathbb{Z} + \frac{j}{2}} V_n$,
and equipped with a linear map
\begin{equation}
V \longrightarrow (\mbox{End}\,V)[[x,x^{-1}]], \qquad v \mapsto Y(v,x)= \sum_{n\in\mathbb{Z}}v_nx^{-n-1},
\end{equation}
and with two distinguished vectors ${\bf 1}\in V_0$, (the {\it vacuum vector})
and $\omega\in V_2$
(the {\it conformal element}) satisfying the following conditions for $u, v \in
V$:  $u_nv=0$ for $n$ sufficiently large; $Y({\bf 1},x)v=v$; $Y(v,x){\bf 1}\in V[[x]]$, and $\lim_{x\to
0}Y(v,x){\bf 1}=v$;
\begin{eqnarray}
x^{-1}_0\delta\left(\frac{x_1-x_2}{x_0}\right) Y(u,x_1)Y(v,x_2) - (-1)^{|u||v|}
x^{-1}_0\delta\left(\frac{x_2-x_1}{-x_0}\right) Y(v,x_2)Y(u,x_1) \ \ \ \ \ \nonumber \\
= x_2^{-1}\delta \left(\frac{x_1-x_0}{x_2}\right) Y(Y(u,x_0)v,x_2) \ \ \ \ \ 
\end{eqnarray}
(the {\it Jacobi identity}), where $|v| = j$ if $v \in V^{(j)}$ for $j \in \mathbb{Z}_2$; writing $Y(\omega,x)=\sum_{n\in\mathbb{Z}}L(n)x^{-n-2}$, i.e., $L(n)=\omega_{n+1}$, for $n\in \mathbb{Z}$, then the $L(n)$ give a representation of the Virasoro algebra with central charge $c \in \mathbb{C}$ (the {\it central charge} of $V$); for $n \in  \frac{1}{2} \mathbb{Z}$ and $v\in V_n$, then $L(0)v=nv=(\mbox{wt}\,v)v$; and the {\it $L(-1)$-derivative property} holds: 
$\frac{d}{dx}Y(v,x)=Y(L(-1)v,x)$.
\end{definition}

An {\it automorphism} of a VOSA, $V$, is a linear map $g$ from $V$ to itself, preserving ${\bf 1}$ and $\omega$ such that
the actions of $g$ and $Y(v,x)$ on $V$ are compatible in the sense
that $g Y(v,x) g^{-1}=Y(gv,x)$,
for $v\in V.$ Then $g V_n\subset V_n$ for $n\in  \frac{1}{2} \mathbb{Z}$. 

Let $\mathbb{Z}_+$ denote the positive integers.  If $g$
has finite order, $V$ is a direct sum of the eigenspaces $V^j$ of $g$, i.e., 
$V=\bigoplus_{j\in \mathbb{Z} /k \mathbb{Z} }V^j$,
where $k \in \mathbb{Z}_+$ is a period of $g$ (i.e., $g^k = 1$) and
$V^j=\{v\in V \; | \; g v= \eta^j v\}$,
for $\eta$ a fixed primitive $k$-th root of unity.

\begin{definition} Let $(V,Y,{\bf 1},\omega)$ be a VOSA
and $g$ an automorphism of $V$ of period $k \in
\mathbb{Z}_+$. A {\it weak $g$-twisted $V$-module} is a 
vector space $M$ equipped with a linear map
\begin{equation}
V \longrightarrow  (\mbox{End}\,M)[[x^{1/k},x^{-1/k}]], \qquad 
v \mapsto  Y^g (v,x)=\sum_{n\in \frac{1}{k}\mathbb{Z} }v_n^g x^{-n-1},
\end{equation}
with $v_n^g \in (\mathrm{End} \, M)^{(|v|)}$, and satisfying the following conditions for $u,v\in V$ and $w\in M$: $v_n^g w=0$ for $n$ sufficiently large; $Y^g ({\bf 1},x)w=w$; 
\begin{eqnarray}\label{twisted-jacobi}
x^{-1}_0\delta\left(\frac{x_1-x_2}{x_0}\right)
Y^g(u,x_1)Y^g(v,x_2)-  (-1)^{|u||v|} x^{-1}_0\delta\left(\frac{x_2-x_1}{-x_0}\right) Y^g
(v,x_2)Y^g (u,x_1) \ \ \nonumber \\
= x_2^{-1}\frac{1}{k}\sum_{j\in \mathbb{Z} /k \mathbb{Z}}
\delta\left(\eta^j\frac{(x_1-x_0)^{1/k}}{x_2^{1/k}}\right)Y^g (Y(g^j
u,x_0)v,x_2) \ \ \ \ \ \  \ \ \ \ \ \ \
\end{eqnarray}
(the {\it twisted Jacobi identity}) where $\eta$ is a fixed primitive
$k$-th root of unity.  
\end{definition}

As a consequence of the definition, we have that 
$Y^g(v,x)=\sum_{n\in  \mathbb{Z} + \frac{j}{k}}v_n^g x^{-n-1}$ for $j\in 
\mathbb{Z}/k\mathbb{Z}$ and $v\in V^j$, and for $v \in V$,
we have
$Y_g(gv,x) = \lim_{x^{1/k} \to \eta^{-1} x^{1/k}} Y_g(v,x)$.
It also follows that writing $Y^g (\omega,x)=\sum_{n\in \mathbb{Z}} L^g (n)x^{-n-2}$, i.e., setting $L^g (n)=\omega_{n+1}^g$, for $n\in \mathbb{Z}$, then the $L^g(n)$ satisfy the relations for the Virasoro algebra with central charge $c$ the central charge of $V$.

If we take $g=1$, then we obtain the notion of weak $V$-module.  The term ``weak" means we are making no assumptions about a grading on $M$.  

A {\em weak admissible}  $g$-twisted $V$-module is a weak $g$-twisted
$V$-module $M$ which carries a $\frac{1}{2k}\mathbb{N}$-grading $
M=\bigoplus_{n\in\frac{1}{2k}\mathbb{N}}M(n)$,
such that $v^g_mM(n)\subseteq M(n+\mathrm{wt} \; v-m-1)$ for homogeneous $v\in V$, $n \in \frac{1}{2k} \mathbb{N}$, and $m \in \frac{1}{k} \mathbb{Z}$.  
If $g = 1$, then a weak admissible $g$-twisted $V$-module is called a weak admissible $V$-module.  

An (ordinary) $g$-twisted $V$-module is a weak $g$-twisted $V$-module $M$ graded by $\mathbb{C}$ induced by the spectrum of $L(0).$ That is, we have
$M=\bigoplus_{\lambda \in{\mathbb{C}}}M_{\lambda}$,
where $M_{\lambda}=\{w\in M|L(0)^gw=\lambda w\}$, for $L(0)^g = \omega_1^g$. Moreover we require that $\dim
M_{\lambda}$ is finite and $M_{n/2k +\lambda}=0$ for fixed $\lambda$ and for all
sufficiently small integers $n$.  If $g = 1$, then a $g$-twisted $V$-module is a $V$-module.

\section{The construction and classification of $(1 \; 2 \; \cdots \; k)$-twisted $V^{\otimes k}$-modules}

Let $V=(V,Y,{\bf 1},\omega)$ be a VOSA, and let $k$ be a
fixed positive integer.  Then $V^{\otimes k}$ is also a VOSA, and the permutation group $S_k$
acts naturally on $V^{\otimes k}$ as signed automorphisms, i.e.,  as a right action we have
\begin{eqnarray}
 (1 \; 2 \;  \cdots \;  k) : V \otimes V \otimes \cdots \otimes V  &\longrightarrow &  V \otimes V \otimes \cdots \otimes V\\
v_1 \otimes v_2 \otimes \cdots \otimes v_k  & \mapsto &   (-1)^{|v_1|(|v_2| + \cdots + |v_k|)} v_2 \otimes v_3 \otimes \cdots \otimes v_k \otimes v_1. \nonumber
\end{eqnarray}

Let $g=(1 \; 2 \; \cdots \; k)$.  Below,  we will recall the classification and construction of $g$-twisted $V^{\otimes k}$-modules from \cite{B-odd-superpermutation} and \cite{BV-even-superpermutation}.  This construction is based on a certain operator $\Delta_k(x)$ first defined in \cite{BDM}, which we now recall.

Consider the polynomial $\frac{1}{k} (1 + x)^k - \frac{1}{k}$ in $x \mathbb{Q} [x]$. 
Following \cite{BDM}, for $k \in \mathbb{Z}_+$, we define $a_j \in \mathbb{Q}$ for $j \in \mathbb{Z}_+$, by  
\begin{equation}\label{define-a}
\exp \Biggl( - \sum_{j \in \mathbb{Z}_+} a_j  x^{j + 1} \frac{\partial}{\partial x}
\Biggr) \cdot x = \frac{1}{k} (1 + x)^k -
\frac{1}{k} .
\end{equation}
For example, $a_1=(1-k)/2$ and $a_2=(k^2-1)/12.$
Let $V= (V, Y, {\bf 1}, \omega)$ be a VOSA.  In 
$(\mathrm{End} \;V)[[x^{1/2k}, x^{-1/2k}]]$, define 
\begin{equation}\label{Delta-for-a-module}
\Delta_k (x) = \exp \Biggl( \sum_{j \in \mathbb{Z}_+} a_j x^{- \frac{j}{k}} L(j) 
\Biggr) (k^\frac{1}{2})^{-2L(0)} \left(x^{\frac{1}{2k}\left( k-1\right)}\right)^{- 2L(0)} .
\end{equation}

For $v\in V$, and $k$ any positive integer, denote by $v^j\in V^{\otimes k}$, for $j = 1,\dots, k$, the vector whose $j$-th tensor factor is $v$ and whose other tensor factors are ${\bf 1}$.  Then for $g = (1 \; 2 \; \cdots \; k)$, we have $gv^j=v^{j-1}$ for $j=1,\dots,k$ where $0$ is understood to be $k$. 

Let $(M, Y_M)$ be a $V$-module, and $(M_\sigma, Y_\sigma)$ a $\sigma$-twisted $V$-module, where $\sigma$ is the parity map on $V$.  We define the $g$-twisted vertex operators for $V^{\otimes k}$ on $M$, for $k$ odd, and on $M_\sigma$, for $k$ even, as follows:
Set 
\begin{equation}\label{define-twisted-operators1}
Y_g(v^1, x) = \left\{ \begin{array}{ll}
Y_M(\Delta_k(x)v, x^{1/k}) & \quad \mbox{for $k$ odd}\\
\\
Y_\sigma(\Delta_k(x)v, x^{1/k}) & \quad \mbox{for $k$ even}
\end{array} \right.
\end{equation}
and for $j = 0, \dots, k-1$, define
\begin{equation}\label{define-twisted-operators2}
Y_g(v^{j+1}, x) = \lim_{x^{1/k} \to \eta^j x^{1/k}} Y_g(v^1, x).
\end{equation}

Let $V$ be an arbitrary VOSA and $h$ an automorphism of $V$ of finite order.  Denote the categories of weak, weak admissible and ordinary $h$-twisted $V$-modules by $\mathcal{ C}^h_w(V),$ $\mathcal{ C}^h_a(V)$ and $\mathcal{ C}^h(V)$, respectively.  If $h=1$, we habitually remove the index $h.$

Now again consider the VOSA, $V^{\otimes k}$, and the
$k$-cycle $g = (1 \; 2 \; \cdots \; k)$.  For $k$ odd, define
\begin{equation}
T_g^k \ : \ \mathcal{ C}_w(V) \ \longrightarrow \ \mathcal{ C}^g_w(V^{\otimes k}), \qquad 
  (M,Y_M) \ \mapsto \ (T_g^k(M),Y_g) = (M,Y_g).
\end{equation}
For $k$ even, define
\begin{equation}\label{even-functor}
T_g^k \ : \ \mathcal{ C}^\sigma_w(V) \ \longrightarrow \ \mathcal{ C}^g_w(V^{\otimes k}), \qquad 
  (M_\sigma,Y_\sigma)  \ \mapsto \ (T_g^k(M_\sigma),Y_g) = (M_\sigma,Y_g).
\end{equation}

The following theorem is proved in \cite{B-odd-superpermutation} for $k$ odd, and in \cite{BV-even-superpermutation} for $k$ even.

\begin{theorem}\label{zeroth-theorem} (\cite{B-odd-superpermutation}, \cite{BV-even-superpermutation}) 

(1)  For $k$ odd, the functor $T_g^k$ is an isomorphism {}from the category $\mathcal{ C}_w(V)$ of weak $V$-modules to the category $\mathcal{ C}^g_w(V^{\otimes k})$ of weak $g = (1 \; 2 \; \cdots \; k)$-twisted $V^{\otimes k}$-modules.

(2)  For $k$ even, the functor $T_g^k$ is an isomorphism {}from the category $\mathcal{ C}^\sigma_w(V)$ of weak parity-twisted $V$-modules to the category $\mathcal{ C}^g_w(V^{\otimes k})$ of weak $g = (1 \; 2 \; \cdots \; k)$-twisted $V^{\otimes k}$-modules. 

(3) For any $k \in \mathbb{Z}_+$, the functor $T_g^k$ restricted to the respective subcategories of weak admissible, ordinary or irreducible modules in $\mathcal{C}_w(V)$ or $\mathcal{C}_w^\sigma(V)$, respectively, is an isomorphism between these subcategories and the corresponding subcategory of weak admissible, ordinary or irreducible $g$-twisted $V^{\otimes k}$-modules.  
\end{theorem}

\section{N=2 supersymmetric VOSAs, Ramond sectors, and mirror-twisted sectors}

In this section, we recall the notions of N=1 or N=2 supersymmetric VOSA,  following the notation and terminology of, for instance, \cite{B-vosas}, \cite{B-iso-thm} and \cite{B-n2axiomatic}.  First we will need the notion of several superextensions of the Virasoro algebra.  

The {\it N=1 Neveu-Schwarz algebra} or {\it N=1 superconformal algebra} is the Lie superalgebra with basis consisting of the central element $d$, even elements $L_n$ for $n \in \mathbb{Z}$, and odd elements $G_{r}$ for $r \in \mathbb{Z} + \frac{1}{2}$, and supercommutation relations  
\begin{eqnarray}
\left[L_m ,L_n \right] &=& (m - n)L_{m + n} + \frac{1}{12} (m^3 - m) \delta_{m + n 
, 0} \; d , \label{Virasoro-relation-N1} \\
\left[ L_m, G_{r} \right] &=& \left(\frac{m}{2} - r\right) G_{m+r} , \qquad  \ \left[ G_{r} , G_{s} \right] \, =\, 2L_{r + s}  + \frac{1}{3} \left(r^2 - \frac{1}{4} \right) \delta_{r+s , 0} \; d ,  \ \ \label{N1-Neveu-Schwarz-relation-last}
\end{eqnarray}
for $m, n \in \mathbb{Z}$, and $r,s \in \mathbb{Z} + \frac{1}{2}$.   The {\it N=1 Ramond algebra} is the Lie superalgebra with basis the central element $d$, even elements $L_n$ for $n \in \mathbb{Z}$,  and odd elements $G_r$ for $r \in \mathbb{Z}$, and supercommutation relations given by (\ref{Virasoro-relation-N1})--(\ref{N1-Neveu-Schwarz-relation-last}), where now $r, s \in \mathbb{Z}$.

The {\it N=2 Neveu-Schwarz Lie superalgebra} or {\it N=2 superconformal algebra} is the Lie superalgebra with basis consisting of the central element $d$, even elements $L_n$ and $J_n$ for $n \in \mathbb{Z}$, and odd elements $G^{(j)}_r$ for $j = 1,2$ and $r \in \mathbb{Z} + \frac{1}{2}$, and such that the supercommutation relations are given as follows:  $L_n$, $d$ and $G^{(j)}_r$ satisfy the supercommutation relations for the N=1 Neveu-Schwarz algebra (\ref{Virasoro-relation-N1})--(\ref{N1-Neveu-Schwarz-relation-last}) for both $G_r = G^{(1)}_r$ and for $G_r = G^{(2)}_r$;  the remaining relations are given by
\begin{eqnarray}
\left[L_m, J_n \right] &=& -n J_{m+n}, \qquad \ \left[ J_m, J_n \right] \ =  \ \frac{1}{3} m \delta_{m+n,0} d \label{N2-evens}\\
\left[ J_m, G^{(1)}_{r} \right]  &=&   - i G^{(2)}_{m+r} , \quad  \left[ J_m, G^{(2)}_{r} \right]  \,  = \,  i G^{(1)}_{m+r} ,
\quad \left[ G^{(1)}_r, G^{(2)}_s \right] 
\, =  \,   i  (s-r) J_{r+s} . \ \ \ \  \ \ \label{nonhomo-J-relation} 
\end{eqnarray}
The {\it N=2 Ramond algebra} is the Lie superalgebra with basis consisting of the central 
element $d$, even elements $L_n$ and $J_n$ for $n \in \mathbb{Z}$, and odd elements 
$G_r^{(j)}$ for $r \in \mathbb{Z}$ and $j=1,2$, and supercommutation relations given by those 
of the N=2 Neveu-Schwarz algebra but with $r,s \in \mathbb{Z}$, instead of $r,s \in 
\mathbb{Z} + \frac{1}{2}$.

Note that there is an automorphism of the N=2 Neveu-Schwarz algebra given by
\begin{equation}\label{mirror-map-nonhomo}
\kappa: \ \   G^{(1)}_r \mapsto G^{(1)}_r, \quad G^{(2)}_r \mapsto - G^{(2)}_r,  \quad J_n \mapsto - J_n, \quad  L_n \mapsto L_n,  \quad d \mapsto d,
\end{equation}
called the {\it mirror map} automorphism of the N=2 Neveu-Schwarz algebra.

Let $(V, Y, {\bf 1}, \omega)$ be a VOSA, and suppose there exists $\tau \in V_{3/2}$ such that writing
$Y(\tau, z) = \sum_{n \in \mathbb{Z}} \tau_n x^{-n-1} = \sum_{n \in \mathbb{Z}} G(n + 1/2) x^{-n-2}$,
the $G(n + 1/2) = \tau_{n+1} \in (\mathrm{End} \, V)^{(1)}$ generate a representation of the 
N=1 Neveu-Schwarz Lie superalgebra such that the $L(n)$ are the modes of $\omega$.
Then we call $(V, Y, {\bf 1}, \tau)$ an {\it N=1 Neveu-Schwarz VOSA}, or an {\it N=1 supersymmetric VOSA}, or just an N=1 VOSA for short.

Suppose a VOSA, $V$, has two vectors $\tau^{(1)}$ and $\tau^{(2)}$ such that $(V, Y, \mathbf{1}, \tau^{(j)})$ is an N=1 VOSA for both $j=1$ and $j = 2$, and the $\tau^{(j)}_{n+1} = G^{(j)}(n + 1/2)$ generate a representation of the N=2 Neveu-Schwarz Lie superalgebra.  Then we call such a VOSA an {\it N=2 Neveu-Schwarz VOSA} or an {\it N=2 supersymmetric VOSA}, or for short, an N=2 VOSA.   



For the case of the parity map, $\sigma$, a $\sigma$-twisted $V$-module, for $V$ an N=1 or N=2 VOSA, is naturally a representation of the N=1 or N=2 Ramond algebra, respectively. (See for instance  \cite{B-n2twisted}, \cite{B-varna2011}, as well as references therein). 

Suppose $V$ is an N=2 VOSA such that $V$ has an automorphism $g_\kappa$ which is a lift of the mirror map $\kappa$ for the N=2 Neveu-Schwarz algebra.  That is letting $g_\kappa$ act by conjugation on $\mathrm{End} \; V$, then $g_\kappa$ restricts to the mirror map $\kappa$ on the elements $L(n)$, $J(n)$, and $G^{(j)}(r)$, for $n \in \mathbb{Z}$, $j = 1, 2$, and $r \in \mathbb{Z} + \frac{1}{2}$, which give the N=2 Neveu-Schwarz algebra representation on the N=2 VOSA, $V$.    Following \cite{B-n2twisted}, \cite{B-varna2011}, we call such an automorphism $g_\kappa$ of an N=2 VOSA, $V$, a {\it mirror map}.  Then a $g_\kappa$-twisted $V$-module is naturally a representation of the ``mirror-twisted N=2 Neveu-Schwarz algebra".  The {\it mirror-twisted N=2 Neveu-Schwarz algebra} is the Lie superalgebra with basis consisting of even elements $L_n$, and $J_r$ and central element $d$, odd elements $G^{(1)}_r$ and $G^{(2)}_n$, for $n \in \mathbb{Z}$ and $r \in \mathbb{Z} + \frac{1}{2}$, and supercommutation relations given as follows:  The $L_n$ and $G^{(1)}_r$ satisfy the supercommutation relations for the N=1 Neveu-Schwarz algebra with central charge $d$;  the $L_n$ and $G^{(2)}_n$ satisfy the supercommutation relations for the N=1 Ramond algebra with central charge $d$;  and the remaining supercommutation relations are 
\begin{equation}\label{mirror-twisted-NS}
\left [L_n, J_r\right]  \ = \  -r J_{n+r}, \quad   \left[ J_r, J_s \right]   =   \frac{1}{3} r \delta_{r+s,0} d, \quad  \left[ G^{(1)}_r, G^{(2)}_n \right]   \  = \  - i  (r-n) J_{r+n} 
\end{equation}
\begin{equation}
\left[ J_r, G^{(1)}_s \right] \ = \ - i G^{(2)}_{r+s} ,  \qquad  \left[ J_r, G^{(2)}_n \right]   \   =   \  i G^{(1)}_{r+n} . \label{mirror-twisted-NS-last}
\end{equation}
Note that this mirror-twisted N=2 Neveu-Schwarz algebra is not isomorphic to the ordinary N=2 Neveu-Schwarz algebra \cite{SS}.

\section{Mirror-twisted modules for the class of N=2 VOSAs of the form $V \otimes V$}

There are large classes of N=2 VOSAs of the form $V \otimes V$ such that $V$ is an N=1 VOSA, and $\tilde{\kappa} = (1 \; 2)$, the signed transposition map given by Eqn. (\ref{mirror-map}), is a mirror map for $V \otimes V$.  Examples of such N=2 VOSAs, were studied in \cite{B-n2twisted}.  These include the following examples:  Let $V_L$ be a rank $d$ positive definite integral lattice VOSA or the $d$ free boson vertex operator algebra, and let $V_{fer}^d$ be the $d$ free fermion VOSA.  As noted in \cite{B-n2twisted}, the VOSA $V = V_L \otimes V_{fer}^d$, is naturally an N=1 VOSA, and $V \otimes V$ is naturally an N=2 VOSA.  This uses the construction of a VOSA from a positive definite integral lattice, following for instance \cite{DL}, \cite{X}, \cite{BV-fermions}.    Such N=2 VOSAs have more than one mirror map as was shown in \cite{B-n2twisted}, where the author constructed mirror-twisted modules for these VOSAs for the other mirror map.  

For such N=2 VOSAs of the form $V \otimes V$, and for the signed transposition mirror-map $\tilde{\kappa}$, we have the following immediate corollary to Theorem 1.  

\begin{corollary}  The category of weak mirror-twisted $(V\otimes V)$-modules for the signed transposition mirror map automorphism of an N=2 VOSA of the form $V \otimes V$ is isomorphic to the category of weak N=1 Ramond-twisted $V$-modules (i.e., parity-twisted $V$-modules).  In addition, the subcategories of weak admissible, ordinary, or irreducible modules are isomorphic.   
\end{corollary}

In particular, it follows that if $M_\sigma$ is a representation of the N=1 Ramond algebra such that $M_\sigma$ is a weak parity-twisted module for an N=1 VOSA, $V$, and such that $V \otimes V$ is an N=2 VOSA, then $M_\sigma$ is also naturally a representation of the mirror-twisted N=2 superconformal algebra and is a weak $\tilde{\kappa}$-twisted module for $V \otimes V$.   

Furthermore, from the construction of such modules given by the functor $T_g^k$ for $k = 2$ as in (\ref{define-twisted-operators1}), (\ref{define-twisted-operators2}), (\ref{even-functor}), (see also \cite{BV-even-superpermutation}), we have as a consequence of Corollary 6.5 in \cite{BV-even-superpermutation}, the following:

\begin{corollary}\label{graded-dimension-corollary} $M_{\tilde{\kappa}}$ is an ordinary $\tilde{\kappa}$-twisted $(V \otimes V)$-module with graded dimension 
\[\mathrm{dim}_q M_{\tilde{\kappa}} =tr_{M_{\tilde{\kappa}}} q^{-2c/24 + L^{\tilde{\kappa}}(0)} = q^{-c/12} \sum_{\lambda \in \mathbb{C}} \mathrm{dim} (M_\lambda) q^\lambda \] 
if and only if $(T_{\tilde{\kappa}}^2)^{-1} (M_{\tilde{\kappa}}) = M_{\tilde{\kappa}}$ is an ordinary $\sigma$-twisted $V$-module with graded dimension
\[ \mathrm{dim}_q (T_{\tilde{\kappa}}^2)^{-1} (M_{\tilde{\kappa}}) =  tr_{M_{\tilde{\kappa}}}  q^{-c/24 + L^\sigma (0)} =  \mathrm{dim}_{q^2} M_{\tilde{\kappa}}, \]
where $c$ is the central charge of $V$.
\end{corollary}

\begin{acknowledgement}
The author would like to thank Vladimir Dobrev and the other organizers for their kind invitation to present her research at the workshop in Varna, and for their generous hospitality during her stay.
\end{acknowledgement}



\begin{thebibliography}{BDM..}

%

\bibitem[AM]{Milas-Adamovic2008} D. Adamovi\'c and A. Milas, 
The $N=1$ triplet vertex operator superalgebras: twisted sector, {\it 
SIGMA Symmetry Integrability Geom. Methods Appl.} {\bf 4} (2008), Paper 087, 24 pp. 

\bibitem[B1]{B-vosas} K. Barron, ``$N=1$ Neveu-Schwarz vertex operator superalgebras over Grassmann algebras and with odd formal variables" in {\it Representations and Quantizations: Proceedings of the International Conference on Representation Theory, July 1998}, Shanghai China, Springer-Verlag, 9--39.

%
\bibitem[B2]{B-iso-thm} K. Barron, The notion of $N=1$ supergeometric vertex operator superalgebra and the isomorphism theorem, {\it Commun. in Contemp. Math.},  Vol. 5, No. 4, (2003), 481--567.


\bibitem[B3]{B-n2axiomatic} K. Barron, Axiomatic aspects of N=2 vertex superalgebras with odd formal variables, {\it Commun. in Alg.} {\bf 38} (2010), 1199--1268.




\bibitem[B4]{B-n2twisted} K. Barron, Twisted modules for N=2 supersymmetric vertex operator superalgebras arising from finite automorphisms of the N=2 Neveu-Schwarz algebra, arXiv:1110.0229.

\bibitem[B5]{B-varna2011} K. Barron, On twisted modules for N=2 supersymmetric vertex operator superalgebras, in:  ``Proceedings of the IXth International Workshop on Lie Theory and Its Applications in Physics", June 2011, Varna, Bulgaria; ed. V. Dobrev, Springer 2013, 411--420.

\bibitem[B6]{B-odd-superpermutation} K. Barron, Twisted modules for tensor product vertex operator superalgebras and permutation automorphisms of odd order, arXiv:1310.1956.

\bibitem[BDM]{BDM} K. Barron, C. Dong, and G. Mason, Twisted sectors for tensor products vertex operator algebras associated to permutation groups, {\it Comm. Math. Phys.} {\bf 227} (2002), 349--384.

\bibitem[BHL]{BHL} K. Barron, Y.-Z. Huang, and J. Lepowsky, An equivalence of two constructions of permutation-twisted modules for lattice vertex operator algebras, {\it J. Pure Appl. Algebra} {\bf 210} (2007), 797--826.

\bibitem[BV1]{BV-fermions}  K. Barron and N. Vander Werf, On permutation-twisted free fermions and two conjectures, in: ``Proceedings of the XXIst International Conference on Integrable Systems and Quantum Symmetries", June 2013, Prague, Czech Republic; ed. C. Burdik, O. Navratil and S. Posta, {\it Jour. of Physics: Conference Series} {\bf 474} (2013), 012009; 35 pages.

\bibitem[BV2]{BV-even-superpermutation}  K. Barron and N. Vander Werf, Permutation-twisted modules for even order cycles acting on tensor product vertex operator superalgebras, arXiv:1310.3812.

\bibitem[BFK]{BFK1986}  W. Boucher, D. Friedan, and A. Kent, Determinant formulae and unitarity for the N=2 superconformal algebras in two dimensions or exact results on string compactification, {\it Phys. Lett. B} {\bf 172} (1986), no. 3-4, 316--322.

\bibitem[DPZ]{DPZ1985} P. Di Vecchia, J.  Petersen,  and H. Zheng,  N=2 extended superconformal theories in two dimensions, {\it Phys. Lett. B} {\bf 162} (1985), no. 4-6, 327--332. 

\bibitem[Dob]{Dobrev1987} V. Dobrev, Characters of the unitarizable highest weight modules over the N=2 superconformal algebras, {\it Phys. Lett. B} {\bf 186} (1987),  43--51.  See also 	arXiv:0708.1719v1.

\bibitem[DL]{DL} C. Dong and J. Lepowsky, The algebraic structure of relative twisted vertex operators, {\it J. Pure Appl. Algebra} {\bf 110}  (1996), 259--295.

\bibitem[DG2]{DG2001} M. D\"orrzapf, and B. Gato-Rivera, Singular dimensions of the N=2 superconformal algebras. II. The twisted N=2 algebra, {\it Comm. Math. Phys.} {\bf 220} (2001), no. 2, 263--292. 

\bibitem[FFR]{FFR}  A. Feingold, I. Frenkel, and J. Ries, Spinor Construction of Vertex Operator Algebras, Triality, and $E^{(1)}_8$, {\it Contemp. Math.} {\bf 121} (1991).

\bibitem[FQS]{FQS1985}  D. Friedan, Z. Qiu, and S. Shenker, Superconformal invariance in two dimensions and the tricritical ising model, {\it Phys. Lett. B} {\bf 151} (1985), 37--43.

\bibitem[G]{Gato-Rivera2002}  B.  Gato-Rivera,  Construction formulae for singular vectors of the topological and of the Ramond N=2 superconformal algebras,  {\it Internat.  J.  of Modern Phys. A} {\bf 17} (2002), 4515--4541. 

\bibitem[GKO]{GKO1986} P. Goddard, A. Kent, and D. Olive, 
Unitary representations of the Virasoro and super-Virasoro algebras,
{\it Comm. Math. Phys.} {\bf 103} (1986), 105-–119. 

\bibitem[IK1]{IK1999+}  K. Iohara and Y. Koga, R\'esolutions de type Bernstein-Gel'fand-Gel'fand pour les super-alg\`ebres de Virasoro $N=1$, {\it 
C. R. Acad. Sci. Paris S\'er. I Math.} {\bf 328} (1999), 381–-384.  
Representation theory of Neveu-Schwarz and Ramond algebras. I.  Verma modules, {\it Adv. Math.} {\bf  178} (2003),  1–-65. 
Representation theory of Neveu-Schwarz and Ramond algebras. II. Fock modules, {\it Ann. Inst. Fourier (Grenoble)} {\bf 53} (2003), 1755–-1818.
The structure of pre-Verma modules over the $N=1$ Ramond algebra, {\it Lett. Math. Phys.} {\bf 78} (2006), 89–-96. 


Fusion algebras for $N=1$ superconformal field theories through coinvariants. II+$\frac{1}{2}$. Ramond sector, {\it Int. Math. Res. Not. IMRN} {\bf 2009}, 2374-–2416. 

\bibitem[IK2]{IK2004}  K. Iohara and Y. Koga,  Representation theory of $N=2$ super Virasoro algebra: twisted sector, {\it J. of Funct. Anal.} {\bf 214} (2004), 450--518.

\bibitem[K]{Kiritsis1988}  E. Kiritsis,  Character formulae and the structure of the representations of the N=1,N=2 superconformal algebras, {\it Internat. J. Modern Phys. A} {\bf 3} (1988),  1871--1906.

\bibitem[L]{Li-twisted} H. Li, Local systems of twisted vertex operators, vertex operator superalgebras and twisted modules, in: ``Moonshine, the Monster, and related topics (South Hadley, MA, 1994), {\it Contemp. Math} {\bf 193} (1996), 203--236.

\bibitem[LSZ]{LSZ2010}  J. Li, Y. Su, and L. Zhu, Classification of indecomposable modules of the intermediate series over the twisted N=2 superconformal algebra, {\it J. of Math. Phys.} {\bf 51} (2010), no. 8, 083515, 17 pp.

\bibitem[Ma]{Matsuo1987} Y. Matsuo, Character formula of $c<1$ unitary representation of N=2 superconformal algebra, {\it Progr. Theoret. Phys.} {\bf 77} (1987), 793--797.

\bibitem[Mi]{Milas}  A. Milas, Characters, supercharacters and Weber modular functions,  {\it J. Reine Angew. Math.} {\bf 608} (2007), 35–-64.

\bibitem[S]{S} N. Scheithauer,  Vertex algebras, Lie algebras, and superstrings, {\it J. Algebra} {\bf 200} (1998),  363--403.

\bibitem[SS]{SS}  A. Schwimmer and N. Seiberg, Comments on the $N=2,3,4$ superconformal algebras in two dimensions, {\it Phys. Lett. B}, {\bf 184} (1986), 191--196.

\bibitem[X]{X}  X. Xu,  {\it Introduction to vertex operator superalgebras and their modules}, Math. and its Applications, Vol. 456, Kluwer Acad. Pub., Dordrecht, 1998.

\end{thebibliography}
\end{document}